\documentclass[11pt, fleqn]{article}

\usepackage{latexsym}
\usepackage{graphics}
\usepackage{graphicx}
\usepackage{epstopdf}
\usepackage{amsmath}
\usepackage{amsfonts}
\usepackage{eucal}
\usepackage{setspace}

\usepackage[left=1.2in,top=0.9in,right=1in,bottom=1.4in,nohead]{geometry}

\newcommand{\ind}{\setlength{\parindent}{2em}}
\newcommand{\noind}{\setlength{\parindent}{0em}}
\renewcommand{\(}{$ \:}
\renewcommand{\)}{\: $}

\def\h{\hspace{.5em}}

\def\v{\vspace{.5ex}}
\def\vv{\vspace{1ex}}
\def\vvv{\vspace{2ex}}
\def\vvvv{\vspace{4ex}}

\title{\textbf{Exact volume of hyperbolic 2--bridge links}}
\author{\it Anastasiia Tsvietkova}
\date{}

\everymath{\displaystyle}
\begin{document}

\maketitle
 \footnotesize
 \begin{description} \item \qquad \textbf{ Abstract.}  W. Thurston suggested a method for computing hyperbolic volume of hyperbolic 3--manifolds, based on a triangulation of the manifold. The method was implemented by J. Weeks in the program SnapPea, which produces a decimal approximation as a result. For hyperbolic 2-bridge links, we give formulae that allow one to find the exact volume, i.e. to construct a polynomial and to find volume as an analytic function of one of its roots. The computation is performed directly from a reduced, alternating link diagram.

\


\


\end{description}

\normalsize

\vvvv
\noind \textbf{1. Introduction}

\vv
\ind The purpose of this note is to calculate the volume of hyperbolic 2-bridge links exactly, \textit{i.e.} to construct a polynomial such that the volume can be expressed as an analytic function of one of its roots. It turns out that one can do it directly from a reduced, alternating diagram of a 2-bridge link. The constructive process amounts to assigning labels to the link diagram, counting twists and bigons, and substituting the corresponding labels and numbers in the given formulas. This idea emerged from an example, considered by M. Thistlethwaite.

\v

\ind This work was motivated by two questions.  The first one concerns computing the volume of a hyperbolic 3-manifold exactly. W. Thurston suggested a method for computing the volume, based on a triangulation of a manifold (\cite{Thurston}). If done by hand, the process becomes tedious even for manifolds with just several tetrahedra. The method was implemented by J. Weeks in the program SnapPea (\cite{Weeks}), which produces a decimal approximation as a result. The program Snap (\cite{Snap}), intended for exact calculations, followed. It approximates the hyperbolic structure to a high precision, and then makes an intelligent guess of the corresponding algebraic numbers (from which the volume can be computed). For hyperbolic 2-bridge links, we suggest a simple alternative. It does not involve constructing gluing equations for a triangulation (compared to SnapPea) or using the LLL algorithm for guessing the polynomial (compared to Snap). 

\v

The second question concerns relating diagrammatic properties of a link to the geometry of its complement, and, in particular, to its hyperbolic volume. In this spirit, estimates for the volume of various families of links were previously obtained. For example, M. Lackenby showed in \cite{Lackenby} that the hyperbolic volume of an alternating link has upper and lower bounds as functions of the number of twist regions of a reduced, alternating diagram. D. Futer, E. Kalfagianni and J. Purcell extended these results to highly twisted knots (\cite{Futer}) and to sums of alternating tangles (\cite{Futer2009}). We do not provide any bounds (asymptotically sharp volume bounds for 2-bridge links are given in \cite{Bounds}), but the suggested calculation of the exact volume is based solely on the layout of a reduced link diagram.

Previously, polynomials that allow one to compute volumes of 2-bridge links exactly were obtained in \cite{SakumaWeeks} using different methods. In particular, the authors reduced the gluing equations for a triangulation to a single equation. This equation can be obtained by substituting subvectors of $(a_1, a_2, ..., a_n)$, where the numbers $a_1, a_2, ..., a_n$ come from a continued fraction expansion of the link, into a set of recursive formulae given by the authors.  The polynomials suggested here are obtained from a reduced alternating diagram instead, and are closely tied to it. Into the suggested formulae one substitutes labels of the link diagram; these labels are complex numbers assigned to strands and crossings. The construction implies that these polynomials immediately give extra geometric information, such as, for example, translation distance along horoballs between crossing
arcs, or distances and angles between specific horospheres and preimages of meridians in $\mathbb{H}^3$. These geometric parameters are closely related to the invariant trace field of the link, and the connection will be discussed in an upcoming paper (\cite{WithWalter}). For a discussion on how the degrees of the polynomials obtained by the two different methods relate, see Remark 2.4.

\vvvv
 
\noind \textbf{2. Obtaining exact edge and crossing labels}

\vv

\ind In \cite{ThistlethwaiteTsvietkova}, an alternative method for computing the hyperbolic structure of a link complement is described. The key idea of the method is to use isometries of ideal polygons arising from the regions of a link diagram. The method parameterizes the horoball pattern
obtained by lifting cusp neighborhoods to the universal cover $\mathbb{H}^3$ using complex numbers, called edge and crossing labels. According to their names, crossing labels are assigned to crossings of a link projection, and edge labels are assigned to edges of the regions. 

\v
\ind Consider an alternating, reduced diagram $D$ of a hyperbolic 2--bridge link $L$ with $k$ twists, where the leftmost twist has $n_1$ crossings, the twist next to it has $n_2$ crossings, and so on up to $n_k$. We will call the leftmost twist  - the first, the next twist to the right - the second, etc. Note that there always exists a reduced, alternating diagram of a 2-bridge link such that the first and the last twists have at least two crossing each, so we may assume that \(n_1>1\) and $n_k>1$. 

\v
We assume that horospherical cross-sections of the cusps of $S^3\backslash L$ have been chosen
so that a (geodesic) meridian curve on the cross-sectional torus has length 1. The preimage of each cross-sectional torus in the universal cover $\mathbb{H}^3$ is a union of horospheres,
and we specify a complex affine structure on each horosphere by declaring that meridional translation corresponds to the real number 1. Finally, we assume that coordinates in $\mathbb H^{3}$ are chosen so that one of the horospheres is the Euclidean
plane of (Euclidean) height 1 above the $xy$--plane. With these conventions, a crossing label contains the geometric information about a preimage of an arc from an overpass to an underpass of the crossing. An edge label contains information about a preimage of the corresponding arc on the boundary torus. Further details and rigorous definitions of the labels can be found in \cite{ThistlethwaiteTsvietkova}.
\v

Endow $D$ with a labeling scheme  \(\mathcal{L}\) as follows. Every arc from an underpass to an overpass of a crossing is homotopic to any other such arc from the same twist. Therefore, we may use just one crossing label per twist, and edge labels inside any bigon are 0. Denote $w_1$ the leftmost crossing label of the first horizontal twist in $D$ (and hence all crossing labels of that twist), $w_2$ the leftmost crossing label of the second twist (and hence all crossing labels of that twist), and so on up to $w_k$. Recall from \cite{ThistlethwaiteTsvietkova} that given a checkerboard coloring of $D$, if $u$ and $v$ label two different sides of the same edge of $D$, then $u-v=\pm1$ (the sign is plus iff $u$ is in a white region). Hence the edge labels outside bigons are $\pm1$ depending on the coloring chosen. Let us label the rest of the edges of the regions adjacent to twists as shown in Fig. 1. In addition, let us choose an orientation of the link, and let $\kappa_j =1$, if two strands in the $j^\text{{th}}$ twist are oriented coherently, and $\kappa_j=-1$ otherwise.

\v
\begin{center}

\includegraphics[scale=.8]{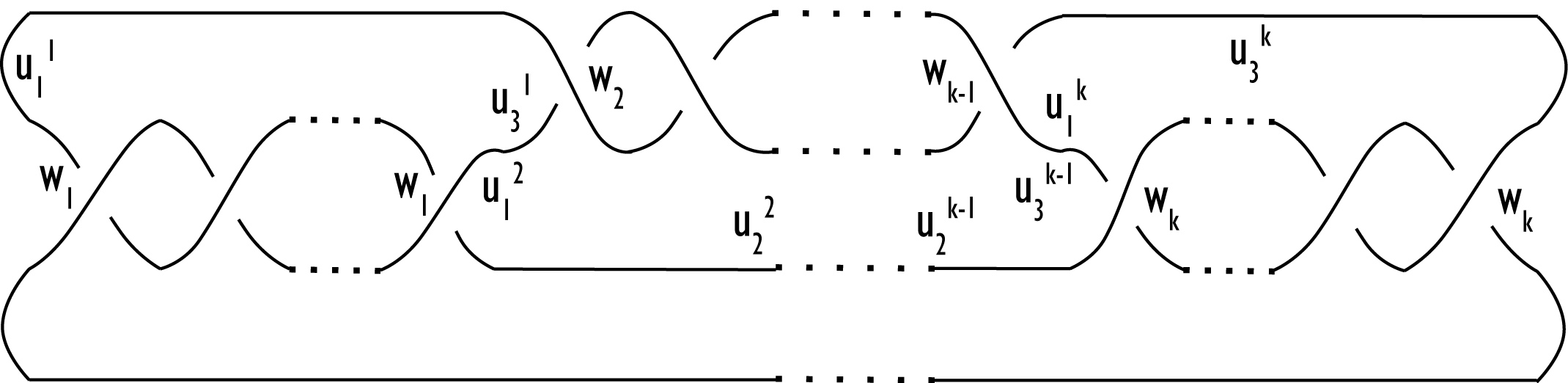}

\v
\it Fig. 1

\end{center}

The layout of a reduced, alternating diagram of a 2-bridge link allows us to find all the labels solely from the label $w_1$. This can be done recursively, using regions of a diagram and proceeding from left to right. The rightmost region yields an extra relation, which gives a polynomial in $w_1$. Therefore, $w_1$ can be found exactly, and all other labels can be expressed in terms of $w_1$. 

\v

The recursive formulae for the labels are a bit unwieldy, but are obtained in a straightforward manner from the region equations given in \cite{ThistlethwaiteTsvietkova}, as we demonstrate in the next proposition. They involve Fibonacci type polynomials defined for the $j^\text{{th}}$ twist as follows:
 $f_0=0$,\h $f_1=1$, \h and \h $f_{m+1}=f_m-\kappa_j w_j
 f_{m-1}$ \h for all natural $m$. Remark 2.5 provides the closed formulae for these polynomials. %The closed formulae for the labels can be obtained as well, though they are not of interest to us in this paper.
 
 \v
 
In this note, we consider edge and crossing labels mainly as parameters aiding in the computation of hyperbolic volume. However the definition of the labels implies that they immediately provide extra geometric information. In particular, \h  $w_i$, \h $i=1, 2, ..., n,$ \h give geometric information about intercusp geodesics at the crossings of $D$, while \h $u_j^i$,\h $i=1, 2, 3,$ \h give information about translations on the boundary torus between the crossings in $D$ (the ones labeled by $w_{j-1}$ and $w_j$, $w_{j-1}$ and $w_{j+1}$, $w_j$ and $w_{j+1}$ respectively). This geometric perspective is explained in \cite{ThistlethwaiteTsvietkova} in detail.

\vv

\noind \textbf{Proposition 2.1.} \h\h Let $D$ be an alternating, reduced diagram  of a hyperbolic 2-bridge link, endowed with the labeling scheme \( \mathcal{L} \) as above. Suppose the $j^\text{{th}}$ twist has \( n \) crossings. Then the labels \( u_3^ j\), \h \( u_1^j\), and \( w_{j+1} \) can be found from \( w_j \) as follows:

\v
\noind \textbf{(i)}   \h\h   for $j=1$,\h $u_1^1=u^1_3=\frac{\kappa_1w_1f_{n-1}}{f_{n}}$ \h and \h $w_2=-\frac{(\kappa_1w_1)^{n}}{\kappa_2f_{n}^2}$;

\v
\noind \textbf{(ii)} \h for \h $1<j<k$  \h and \h $n>1$, \h $u_1^{j}=u_3^{j-1}+(-1)^{j-1}$, \h \(u_3^j=\frac{w_j(\kappa_j u_1^j f_{n-1}(-1)^{j+1}-w_j f_{n-2})}{u_1^j f_{n} + (-1)^{j}\kappa_j w_j f_{n-1}}\), \h  and \h
 $w_{j+1}=\frac{\kappa_{j-1}(\kappa_jw_j)^{n} w_{j-1}}{\kappa_{j+1}(u_1^{j}f_{n}+(-1)^j\kappa_j w_j f_{n-1})^2}$; 
\v

\textbf{(iii)} \h for \h $1<j<k$ \h and \h $n=1$, \h $u_1^{j}=u_3^{j-1}+(-1)^{j-1}$, \h $u_3^j=\frac{\kappa_jw_j}{u_1^j}$,  \newline $w_{j+1}=\frac{\kappa_{j-1}\kappa_jw_{j-1}w_j}{\kappa_{j+1}(u_1^j)^2}$;

\v
\textbf{(iv)} \h for $j=k$, \h $u_3^k=\frac{(-1)^{k+1}\kappa_kw_kf_{n-1}}{f_{n}}$ \h and \h $u^k_1=u_3^{k-1}+(-1)^{k-1}$.

\vv

\noind \textbf{Proof.} \h\h  Let $R_j$ be a region of the diagram $D$, adjacent to the $j^\text{{th}}$ twist with $n$ crossings. The region $R_j$ gives a rise to a disk $\Delta_j$, whose boundary consists of sub-arcs on the peripheral torus
travelling between adjacent crossings incident to $R_j$, and arcs travelling
between the underpass and the overpass at crossings of $R_j$. Denote the corresponding ideal polygon in the cover $\mathbb{H}^3$ by $\tilde{\Delta}_j$. Recall from \cite{ThistlethwaiteTsvietkova} that the shape parameter assigned to the preimage of an arc from an overpass to an underpass is, up to sign, the quotient of the label at that crossing by the product of the two incident edge labels.

\v

\ind Without loss of generality, let us choose such a checkerboard coloring of the regions of the diagram $D$ so that
the bigons in the first twist are black. Recall that for two edge labels $u$ and $v$ located at different sides of the same edge, $u-v=1$ if $u$ is in a white region, and $u-v=-1$ otherwise. This together with the chosen coloring yields the relations for $u_1^j$ from (ii) and (iii), and the relation for $u_1^k$ from (iv). Note also that the edge labels are  \( -1 \) outside a white bigon, and \( 1 \) outside a black bigon, which can be written as $(-1)^{j+1}$ for the $j^\text{{th}}$ twist. 

\v

For every twist, there are two possible situations, namely the number of crossings $n=1$, and \( n > 1\). Suppose first that \(n>1\). Then for $1<j<k$ \h (Fig. 2(i)), some of the shape parameters for $R_j$ are 
\begin{center}

\h $\zeta_1=(-1)^{j+1}\frac{\kappa_jw_j}{u_1^j}$, \h $\zeta_2= \frac{\kappa_{j-1}w_{j-1}}{u_1^j u_2^j}$, \h $\zeta_3=\frac{\kappa_{j+1}w_{j+1}}{u_3^j u_2^j}$, \h and \h $\zeta_4=(-1)^{j+1}\frac{\kappa_jw_j}{u_3^j}$. 

\vvvv

\includegraphics[scale=.6]{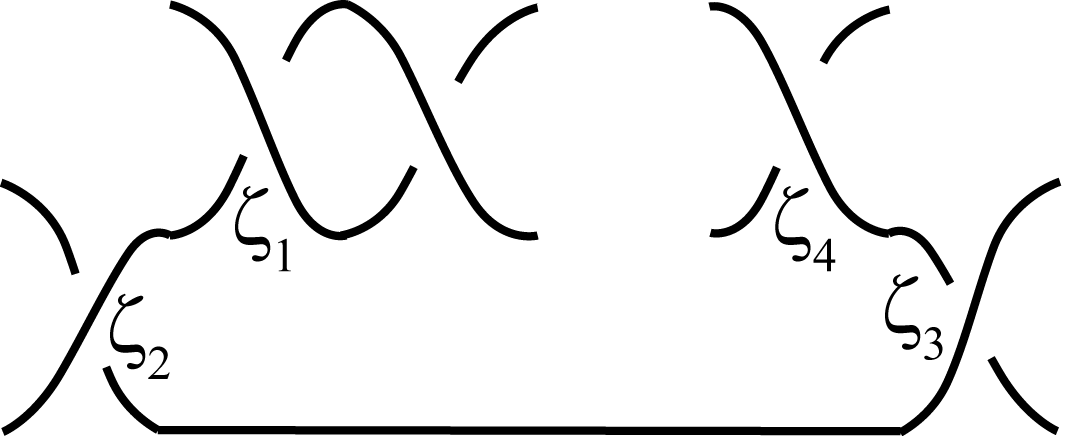} \h\h \h \h \h \h \h
\includegraphics[scale=.55]{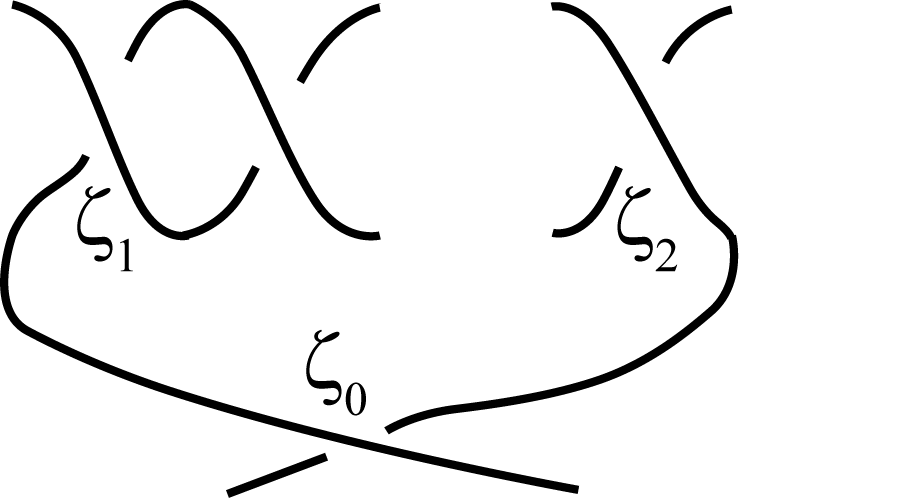}

\end{center}
\textit{\hspace{1.3in}(i) $1<j<k$ \hspace{1.3in} (ii) $j=0$ or $j=k$}
\begin{center}
\it Fig. 2
\end{center}

\noind For the first and last regions the situation is slightly different (Fig. 2(ii)): 
\begin{center} $\zeta_1=\frac{\kappa_1w_1}{u_1^1}$, \h $\zeta_0=\frac{\kappa_2w_2}{u_3^1u_1^1}$,\h and \h $\zeta_2=\frac{\kappa_1w_1}{u_3^1}$ \h for  $R_1$; 

\

$\zeta_1=(-1)^{k+1}\frac{\kappa_kw_k}{u_3^k}$, \h $\zeta_0=\frac{\kappa_{k-1}w_{k-1}}{u_1^k u_3^k}$, \h and \h $\zeta_2=(-1)^{k+1}\frac{\kappa_kw_k}{u_1^k}$ \h for $R_k$.
\end{center}
\noind All the other shape parameters of $R_j$ are  $\kappa_jw_j$  for all $j$. 

\v

\ind Let 
$Z_i = \left[
\begin{array}{cc}
  0 & -\zeta_i \\
  1 & -1 \\
\end{array}
\right]$, \h  \( 0\leq i \leq 4\), \h and \h
 W = $\left[
\begin{array}{cc}
  0 & -\kappa_jw_j \\
  1 & -1 \\
\end{array}
\right]$.

\v

\ind First suppose $1<j<k$, \textit{i.e.} the $j^\text{{th}}$ twist is not the first or the last one.  Denote $c_m$ the element in row 2, column 1 of the product matrix $Z_3 Z_2 Z_1 W^{m-2} Z_4$. Using mathematical induction, one can prove that 
\begin{center}
$c_m=(-1)^{m-1} f_{m}+(-1)^{m} f_{m-1} \zeta_1+(-1)^{m} f_{m} \zeta_2$ \h for all natural \h $m>1$.
\end{center}

\ind The product  \( Z_3 Z_2 Z_1 W^{n-2} Z_4 \) corresponds to the composition of hyperbolic isometries, rotating the polygon $\tilde{\Delta}_j$. Since the polygon closes up, this composition is 1, whence \( Z_3 Z_2 Z_1 W^{n-2} Z_4 \) is a scalar multiple of the identity matrix. Therefore,  $c_n=0$. This together with the above equality for \( c_m \) implies  
 \( \zeta_2=\frac{f_{n}-f_{n-1}\zeta_1}{f_{n}}\). Substituting the shape parameters, we obtain

\begin{center}

\(u_2^j=\frac{\kappa_{j-1} w_{j-1}f_{n}}{u_1^jf_{n}+\kappa_j w_jf_{n-1}(-1)^{j}}\).

\end{center}

\ind Similarly, the product \( Z_2 Z_1 W^{n-2} Z_4 Z_3 \) yields 
\begin{center} $(-1)^{n-1}f_{n}+(-1)^{n-1}\zeta_1\zeta_4f_{n-2}+(-1)^n(\zeta_1+\zeta_4)f_{n-1}=0$, \h whence \h
\(\zeta_4=\frac{f_{n}-f_{n-1}\zeta_1}{f_{n-1}-f_{n-2}\zeta_1}\),
\end{center}

\noind and the equality for \( u_3^j \) stated in (ii) follows (there is a factor $(\kappa_j)^2$ in the resulting formulae that we can safely omit, since $\kappa_j=\pm1$, and therefore $(\kappa_j)^2=1$).

\v

\ind Using the symmetry of the region \( R_j\), we can substitute \( \zeta_1 \) by \( \zeta_4\), and \( \zeta_2 \) by \( \zeta_3 \) in the first equality proved by induction. Then

\begin{center}   
\(\zeta_3=\frac{f_{n}-f_{n-1} \zeta_4}{f_{n}}=1-\frac{f_{n-1}}{f_{n}} \cdot \frac{f_{n}-f_{n-1}\zeta_1}{f_{n-1}-f_{n-2}\zeta_1}\).

\end{center}

\noind Note that 

\v
$f_{m+1}^2-f_{m+2}f_{m}=(f_{m}-\kappa_j w_jf_{m-1})f_{m+1}-(f_{m+1}-\kappa_j w_jf_{m})f_{m}=\kappa_j w_j(f_{m}^2-f_{m+1}f_{m-1})$
\v

 \noind for every natural \( m>1\), and therefore, by induction,  $(\kappa_j w_j)^{m-1}=f_{m}^2-f_{m+1}f_{m-1}$. So the relation for $\zeta_3$ becomes 
$\zeta_3=\frac{(\kappa_jw_j)^{n-2}\zeta_1}{f_{n}(f_{n-1}-f_{n-2}\zeta_1)}$.
 Substituting the shape parameters and the above formulae for \( u_2^j \) and \( u_3^j\), we obtain the relation for \( w_{j+1} \) stated in (ii).

\v

\ind For the first and last regions (\( j=1 \) and \(  j=k \)), the matrix product \( Z_0Z_1W^{n-2}Z_2 \) and the symmetry of the region imply \( \zeta_1=\zeta_2=\frac{f_n}{f_{n-1}}\). The relations (i) and (iv) follow.  

\v

\ind In a similar fashion, the product $W^{n-2}Z_2Z_0Z_1$ yields

\begin{center}
$\zeta_0=\frac{-\zeta_2f_{n-2}+f_{n-1}}{f_{n-1}}=-\frac{f_{n-1}^2-f_{n-2}f_{n}}{f^2_{n-1}}=-\frac{(\kappa_jw_j)^{n-2}}{f^2_{n-1}}$, \h
\end{center}
\noind and the formula for $w_2$ stated in (i) follows.

\v

\ind In the case of  \( n=1\), \( R_j \) is 3--sided, and all the shape parameters are equal to 1. This implies (iii), completing the proof. 
${\Box}$

\vv

\noind \textbf{Remark 2.2.} \h In the proof of Proposition 2.1, the relation \( \zeta_2=\frac{f_{n}}{f_{n-1}} \) for the last region yields an extra equation
\begin{center}

$u_1^kf_{n}+(-1)^{k}\kappa_k w_kf_{n-1}=0$,

\end{center}

\noind where \( n>1 \) is the number of crossings in the last twist. This equation together with the recursive formulae from Proposition 2.1 describe a constructive process of obtaining a polynomial $P$ from the link diagram, for which \( w_1 \) is a root. The polynomial has several complex roots. The root that corresponds to the geometric structure is the one that maximizes hyperbolic volume.

\vv

\ind By a twist link with $n+2$ crossings we will mean a link (or a knot) that has a reduced alternating diagram with two crossings in the first twist and $n$ crossings in the second twist.

\vv
%The root can be used to compute volume of the link complement, as we will demonstrate in the next section.
\noind \textbf{Lemma 2.3.} \h For a twist link with $n+2$ crossings, the degree of the polynomial $P$ is exactly $n$.
\vv

\noind \textbf{Proof.} \h The recursive definition of the polynomial $f_n$ in terms of $w_j$ implies that its degree is 0 for \h $n=0$, \h $\frac{n}{2}-1$ \h for $n$ even, and \h $\frac{n-1}{2}$ \h for $n$ odd.

\ind For a twist link $L$, let $f_n$ be the polynomial for the second twist. From Proposition 2.1, \h $w_2=\pm(w_1)^2$, \h $u_1^1=u_3^1=\pm w_1$, \h $u_1^2=\pm w_1\pm1$, \h and \h $u_3^2=\frac{\pm(w_1)^2f_{n-1}}{f_n}$. \h If we rewrite $f_n$ in terms of $w_1$, its degree will be $n-2$ for $n$ even, and $n-1$ for $n$ odd. Therefore, 
\begin{center}
 \(P(w_1)=(\pm w_1 \pm 1)f_n\pm (w_1)^2f_{n-1}\). 
 \end{center}
 
\noind Treating $n$ even and $n$ odd as two separate cases, we obtain that the degree of $P$ is always $n$. $\Box$ 

\vv 

\noind \textbf{Remark 2.4.} \h Given a link, both the polynomial $P$ from Remark 2.2 and the polynomial obtained using formulae of Sakuma and Weeks (\cite{SakumaWeeks}) can be used for volume computation. The methods and formulae are different, and it is unclear how the two resulting polynomials or their degrees might relate for an arbitrary 2-bridge link. Partial information about the degree of the polynomial suggested by Sakuma and Weeks is given in Lemma II.5.8 of \cite{SakumaWeeks}. In particular, the lemma states that for a twist link with $n+2$ crossings, the degree can be written in terms of a recursive function $\alpha$, defined in \cite{SakumaWeeks}, and is exactly \h $\frac{\alpha(2,n)-1}{2}$ \h if \h $\alpha(2,n)$ is odd. Note that \h $\alpha(2,n)=n\alpha(2)+\alpha(0)=2n+1$ \h is always odd, and therefore the degree of the Sakuma-Weeks polynomial is exactly $n$ for the twist link. This coincides with the degree of the polynomial $P$ (as proved in Lemma 2.3).

\vv 

\noind \textbf{Remark 2.5.} \h The formulae for the labels involve  Fibonacci type polynomials: $f_0=0$,\h $f_1=1$, \h and \h $f_{m+1}=f_m-\kappa w f_{m-1}$ \h for all natural $m$. Below we provide a closed formula for this recurrence relation.

\ind The characteristic polynomial for the recurrence relation is $r^2-r+\kappa w=0$. Solving for $r$, we obtain the characteristic roots \h $r_1=(1+\sqrt{1-4\kappa w})\slash 2$ \h and \h $r_1=(1-\sqrt{1-4\kappa w})\slash 2$. For $r_1$ and $r_2$ to be equal, $w$ would have to be real number, which is never the case for a crossing label $w$ that corresponds to the hyperbolic structure of a link complement. Therefore, $r_1$ and $r_2$ are distinct. Then \h $f_n=C(r_1)^n+D(r_2)^n$, \h where $C$ and $D$ are constants chosen based on two initial relations \h $f_0=0$ \h and \h $f_1=1$. \h Finding $C$ and $D$, one establishes the closed formula
\begin{center}
 $f_n=\frac{(1+\sqrt{1-4\kappa w})^n-(1-\sqrt{1-4\kappa w})^n}{2^n\sqrt{1-4\kappa w}}$.
 \end{center}

\ind This suggests that closed formulae for the labels and the polynomial $P$ might be established as well, albeit they might be rather cumbersome. 

\vvvv

\noind \textbf{3. Shapes of ideal tetrahedra}

\vv

\ind In \cite{SakumaWeeks}, the conjectural canonical cell decomposition of 2-bridge links was described (the proof that this ideal triangulation is indeed the canonical cell decomposition was announced in \cite{Sakuma1} and should appear in a sequel). According to the Sakuma-Weeks description, there are \( 2(n_1+n_2+...+n_k-3) \) tetrahedra occurring
in isometric pairs in the complement \( S^3\backslash L \) of a hyperbolic 2-bridge link $L$. In this section, we give simple formulas that allow one to find the tetrahedra shapes from the labels of a reduced, alternating diagram $D$ of $L$.

\v

Suppose $D$ is endowed with a labeling scheme  \( \mathcal{L} \) as above. We will add more labels, and call the resulting scheme \( \mathcal{L}_1\). On Fig. 3(i), the vertical dotted lines represent geodesics in the complement \( S^3\backslash L \)  that correspond to edges of the canonical cell decomposition. Every \( b_i^j \), except for the first and last ones ( $b_1^1$ and $b_{n_k}^k$ ), is the edge label corresponding to a (directed) Euclidean line
segment on the boundary torus between such a geodesic and the geodesic for the next to the right crossing arc in the twist. The first one, $b_1^1$, corresponds to the segment from the first (leftmost) crossing arc to the second crossing arc. The last one,   $b_{n_k}^k$, corresponds to the segment between the last crossing arc and the one preceding it. In our notation the upper index indicates the adjacent twist. 

\v

Vertical dotted lines subdivide the diagram into parts. We will refer to these parts as ``levels". Examining the Sakuma and Weeks description, we see that there are two isometric tetrahedra on every level.  We will say that an ideal tetrahedron \( T \) is adjacent to a crossing, if an arc from the overpass to the underpass of this crossing is (homotopic to) a truncated edge of \( T\). We will also say that \( T \) is adjacent to a twist, if \( T \) is adjacent to one of the crossings of this twist.
\v

The following theorem gives expressions for the tetrahedra shapes in terms of the labels.

\begin{center}

\includegraphics[scale=.8]{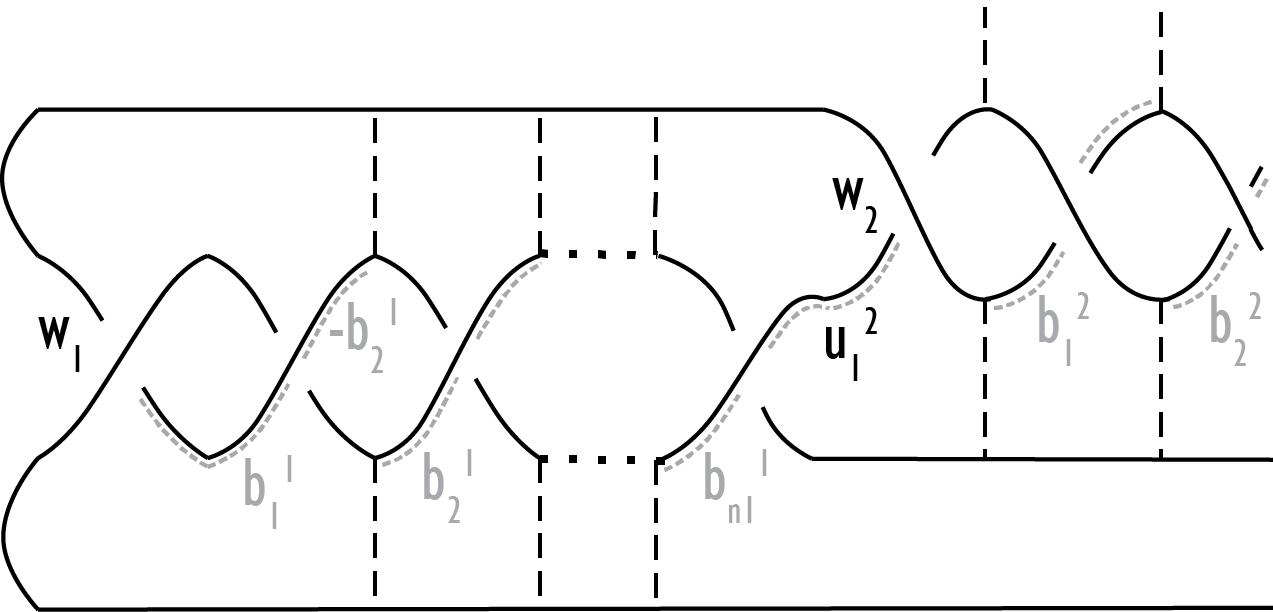}\h\h\h\h\h\h\h\h\h
\includegraphics[scale=.8]{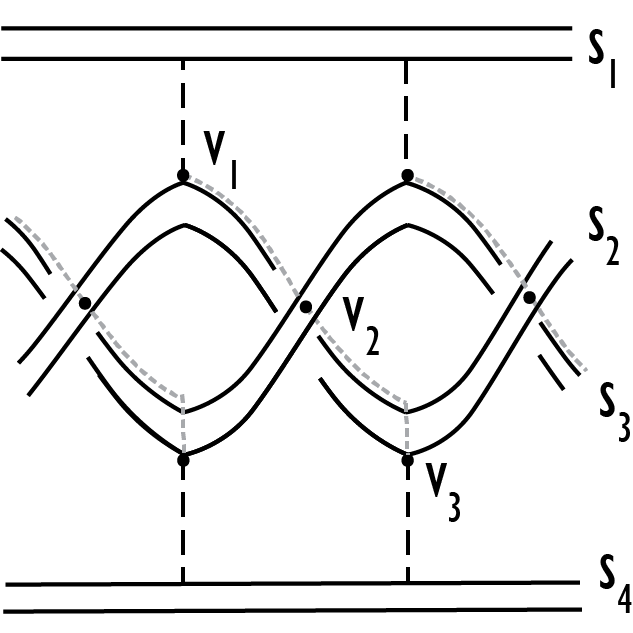}

\end{center}
\textit{\hspace{0.4in}(i) Decomposition of 2--bridge link complement\hspace{0.55in} (ii) Vertices of a cross-section}
\begin{center}
\it Fig. 3
\end{center}

\noind \textbf{Theorem 3.1.} \h  Let \( D \) be the reduced, alternating diagram of a hyperbolic 2-bridge link \( L \) with a labeling scheme \( \mathcal{L}_1 \) as above. In the canonical cell decomposition of \( S^3 \backslash L\), a tetrahedron shape is  \( z \)  if  \( Arg(z)>0\), and  \( 1/z \)  otherwise, where  \( z \)  is a ratio of the form

\v
\noind \textbf{(i)} \h\ \( (b_{n_k-1}^k-1)\slash b_{n_k}^k \) \h for a pair of tetrahedra adjacent to the last crossing of the last twist; 

\v

\noind \textbf{(ii)} \h \( -b_{n_j}^j\slash u_1^{j+1} \) \h\h for a pair of tetrahedra adjacent to the $j^\text{{th}}$ and $(j+1)^\text{{th}}$ twists,  \( n_j>1 \);

\v

\noind \textbf{(iii)} \h \( -u_3^{j}\slash b_1^{j+1} \) \h for a pair of tetrahedra adjacent to the $j^\text{{th}}$ and $(j+1)^\text{{th}}$ twists, \( n_j=1\);

\v

\noind \textbf{(iv)} \h \( b_i^j\slash b_{i+1}^j \) \h\h\h for all other pairs of tetrahedra.

\vv

\noind \textbf{Proof.} \h Denote the strands of the diagram \( D \) by $s_1$, $s_2$, $s_3$, $s_4$. Consider an ideal tetrahedron in the canonical cell decomposition of \( S^3 \backslash L\). Inspection of the Sakuma-Weeks description (\cite{SakumaWeeks}) shows that each Euclidean cusp cross-section on a (thickened) strand $s_i$ has its three vertices \( v_1\), \( v_2\), \( v_3 \) on geodesics joining \( s_i \) with each of the other three
strands. Fig. 3(ii) illustrates this situation for a tetrahedron adjacent to a twist with more than one crossing. Grey dotted arcs from \( v_1 \) to  \( v_2 \) and from  \(v_2 \) to \( v_3 \) along the boundary torus are edges of the cross-section. Fig. 3(i) extends this picture, showing pairs of edges for cross-sections of multiple tetrahedra. 

\v
\ind We have two cases to consider, namely, when a pair of isometric tetrahedra is adjacent just to one twist, and when it is adjacent to two. In the first case, the complex translations corresponding to the edges of a Euclidean cusp cross-section are the edge labels \( b_{i+1} \) and \( -b_i \) (the latter due to the symmetry of the diagram near a twist), unless tetrahedra are adjacent to the last crossing of the last twist. For that last crossing, the complex translations are \( 1-b_{n_{k}-1}^k \) and \( b_{n_k}^k\). In the second case, i.e. at the levels, where the diagram goes from one twist into another, the situation is slightly different. The complex translations along the edges of a Euclidean cross-section correspond to the complex numbers \( u_3^{j} \) and \( b_1^{j+1} \) if the $j^\text{{th}}$ twist has just one crossing (Fig. 4(i)), and to \( b_{n_j}^j \) and  \( u_1^{j+1} \) otherwise (as shown on Fig. 3(i)).

\v

\ind From the layout of the diagram, we see that some crossings are between the strands \( s_1 \) and \( s_2\), while all others are between \( s_2 \) and \( s_3\).  For the first type of crossings, we can write
\( z \) as suggested in (i)-(iv), and then shapes for the corresponding pair of isometric tetrahedra become the usual \h $z$,\h $1-\frac{1}{z}$, \h $\frac{1}{1-z}$. The extra minus before the ratios in (i)-(iv) guarantees that the argument of a shape corresponds to the interior dihedral angle of the tetrahedron (not an exterior one). For the second type of crossing, \( \frac{1}{z} \) gives a corresponding shape. Note that every shape has a positive imaginary part, while its reciprocal has a negative one.
${\Box}$
\vv 

It is intended to express shapes of all tetrahedra in terms of just one label, $w_1$, which is a root of the polynomial $P$ (defined in Remark 2.2). With this purpose in mind, we give the following Lemma.

\begin{center}
\includegraphics[scale=.9]{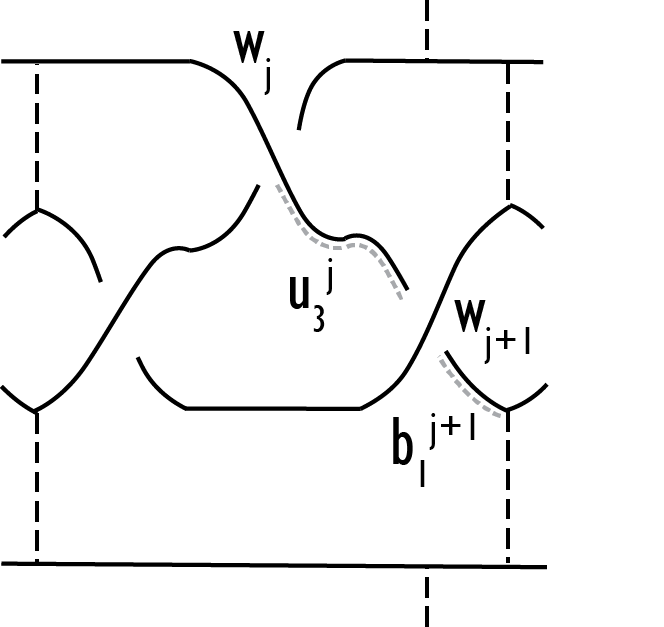} \h
\includegraphics[scale=.9]{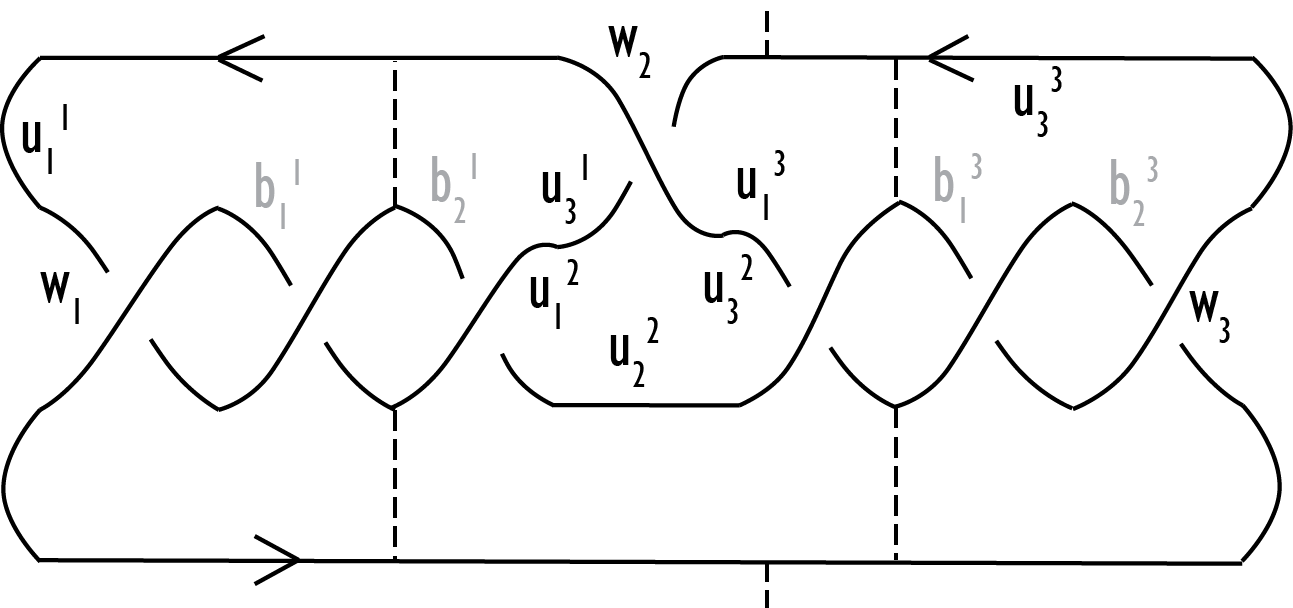}

\end{center}
\vspace{-0.1in}
\textit{(i) Twist with one crossing\hspace{1.5in} (ii) Knot (3, 1, 3)}
\begin{center}
\it Fig. 4
\end{center}

\noind \textbf{Lemma 3.2.} \h Let $D$ be a reduced, alternating diagram of a hyperbolic 2-bridge link, endowed with a labeling scheme  \( \mathcal{L}_1 \) as above. Then every $b_i^j$ can be found as follows: 

\v

\noind \textbf{(i)} \h \h \( b_1^1=1\);

\v

\noind \textbf{(ii)} \h \( b^k_{n_k}=(-1)^{k+1}\);

\v

\noind \textbf{(iii)} \h for all other \( b_j^i\), \( b_{i+1}^j =(-1)^{j+1}-\frac{\kappa_jw_{j}}{b_{i}^j} \) \h and \h   
\( b_1^j=(-1)^{j+1}-\frac{\kappa_jw_j}{u_1^j}$.
\vv

\noind \textbf{Proof.} \h In a region adjacent to the $j^\text{{th}}$ twist, the edge label along each bigon is $(-1)^{j+1}$. Therefore, $b_1^1=1$ \h and \h $b^k_{n_k}=(-1)^{k+1}$.  To find all other $b_i^j$, note that the vertical dotted geodesics split the $j^\text{{th}}$ region into triangles.  All the shape parameters in triangular regions are 1, and equal the quotient of $\kappa_jw_j$ by the product of two incident edge labels. From this we obtain (iii). $\Box$
\vv

\noind \textbf{Corollary 3.1.1.} \h A shape of every tetrahedron in the canonical cell decomposition of $S^3\backslash L$ can be written as a rational function of $w_1$.

\vv
\noind \textbf{Proof.} Theorem 3.1 demonstrates that any tetrahedron shape is a rational function of the labels $w_j$, $u_i^j$, and $b_i^j$. Proposition 2.1 and Lemma 3.2 show that any of these labels, in its turn, can be written as a rational function of $w_1$. $\Box$

\vv

\noind \textbf{Remark 3.3.} \h The dihedral angles of a tetrahedron with the shape \( z_p\), \h \( p \) from 1 to \( 2(n_1+n_2+...+n_k-3)\), are
\( \alpha_p\) = Arg\((z_{p})\), \( \beta_{p}\) = Arg\(\left(1-\frac{1}{z_{p}}\right)\),
   \( \gamma_{p}\) = Arg\(\left(\frac{1}{1-z_{p}}\right)\),
and the volume is thus
\( 2\sum_{p=1}^{2(n_1+n_2+...+n_k-3)}(\Lambda(\alpha_{p}) + \Lambda(\beta_{p})
+ \Lambda(\gamma_{p}))\), where \( \Lambda \) is the Lobachevsky function. Therefore, the volume  of \( S^3\backslash L \) can be calculated solely in terms of \( w_1\).

\vv

\ind We proceed with two examples: one of a knot, and one of a link.

\vv

\noind \textbf{Example 3.4.} \h Consider a 2--bridge knot with a Conway code 3 1 3. There are 4 pairs of isometric tetrahedra in its complement. Fix orientation as on Fig. 4(ii). Note that \h$\kappa_1=\kappa_3=-1$, \h $\kappa_2=1$. 

\v

\ind For the region adjacent to the first (leftmost) twist, \h $f_0=0$, \h $f_1=f_2=1$, \h $f_3=f_2+w_1f_1$. Therefore, from Proposition 2.1, \h $u_1^1=u_3^1=-\frac{w_1f_2}{f_3}$, \h $w_2=\frac{(w_1)^3}{(f_3)^2}$. The region adjacent to the second twist has just one crossing, hence we obtain \h $u_1^2=u_3^1-1$,\h $u_3^2=\frac{w_2}{u_1^2}$,\h $w_3=\frac{w_1w_2}{(u_1^2)^2}$. For the region adjacent to the third twist, $f_0=0$ , \h $f_1=f_2=1$, \h $f_3=f_2+w_3f_1$, and therefore, \h $u_1^3=u_3^2+1$. 

\v

\ind Now we can construct the polynomial in $w_1$ using Remark 2.2 : \h $u_1^3f_3+w_3f_2=0$, \h which becomes \h 
\begin{center}
$1 + 7 w_1 + 18 w_1^2 + 19 w_1^3 + 6 w_1^4 + 2 w_1^5 +  4 w_1^6 - w_1^7=0$.
\end{center}
\noind The root that gives the geometric structure is 
\begin{center}
 $w_1=\frac{-(1196+12\sqrt{177})^{\frac{2}{3}}-112+16(1196+12\sqrt{117})^{\frac{1}{3}}+i\sqrt{3}(1196+12\sqrt{117})^{\frac{2}{3}}-112i\sqrt{3}}{12(1196+12\sqrt{117})^{\frac{1}{3}}}$.

\end{center}
\noind This allows us to compute all other labels exactly as well. 

\v
\ind Let us turn to the edges of the canonical cell decomposition (Fig. 4): \( b_1^1=b_3^2=1\), \h $b_2^1=1+\frac{w_1}{b_1^1}$, \h and \h $b_1^3=1+\frac{w_2}{u_1^3}$ \h by Lemma 3.2. The shape parameters are \( z_1=-\frac{-b_2^1}{b_1^1}\), \h \( z_2=-\frac{u_2^1}{b_2^1}\), \h \( z_3=-\frac{u_3^2}{1-b_1^3}\), \h \( z_4=-\frac{b_3^1-1}{b_2^3}\). If we substitute the decimal approximation of \( w_1\), we obtain the volume of 5.1379412018734177698 .

\vv

\noind \textbf{Example 3.5.} \h Consider a 2--bridge link with a Conway code 3 2 3. Fix orientation as on Fig. 5. For this link, there are 5 pairs of isometric tetrahedra in the canonical cell decomposition. Note that \( \kappa_1=\kappa_3=1\), \( \kappa_2=-1\).

\v

\ind For the first region, \( f_0=0\), \( f_1=f_2=1, \) \( f_3=f_2-w_1f_1\). Therefore, from Proposition 2.1 we obtain \( u_1^1=u_1^3=\frac{w_1 f_2}{f_3}\), \( u_1^2=u_1^1-1\), \( w_2=-\frac{-(w_1)^3}{(f_2)^2}\). For the second region, the Fibonacci-type polynomials \( f_0\), \(f_1\), and \( f_2 \) are the same, and   \( u_1^2=u_3^1-1\), \( u_3^2=\frac{w_2u_1^2f_1-f_0(w_2)^2}{u_1^2f_2-w_2f_1} \), \( w_3=\frac{(w_2)^2w_1}{(u_1^2f_2-w_2f_1)^2}\), \( u_1^3=u_3^2+1\). Using Remark 2.2, we obtain a polynomial:
\begin{center}
$-1+7w_1-18w_1^2+16w_1^3+9w_1^9-19w_1^5-4w_1^6+10w_1^7+4w_1^8=0$.
\end{center}
\noind A decimal approximation of the root that gives the hyperbolic structure is \newline \( w_1=0.45899397977032988781 +
0.2236389499547826586251180*i\). Now we can go back and compute the other labels.

\v

\ind Turn to the edges of the canonical cell decomposition: \(b_1^1=b_2^3=1\), \( b_2^1=1-\frac{w_1}{b_1^1}\), \( b_1^2=-1+\frac{w_2}{u_1^2}\), \( b_1^3=1-\frac{w_3}{u_1^3} \) by Remark 3.2.  The shapes for five pairs of isometric tetrahedra are \( z_1=-(\frac{b_1^1}{-b_2^1})\), \( z_2=-\frac{b_2^1}{u_1^2} \), \(  z_3=-\frac{b_1^2}{u_1^3}\), \( z_4=-\frac{1-b_1^3}{u_3^2} \), \( z_5=-\frac{b_2^3}{b_1^3-1}\). Now we can compute the volume; its decimal approximation is 7.5176898964745685429. 

\begin{center}
\includegraphics[scale=.9]{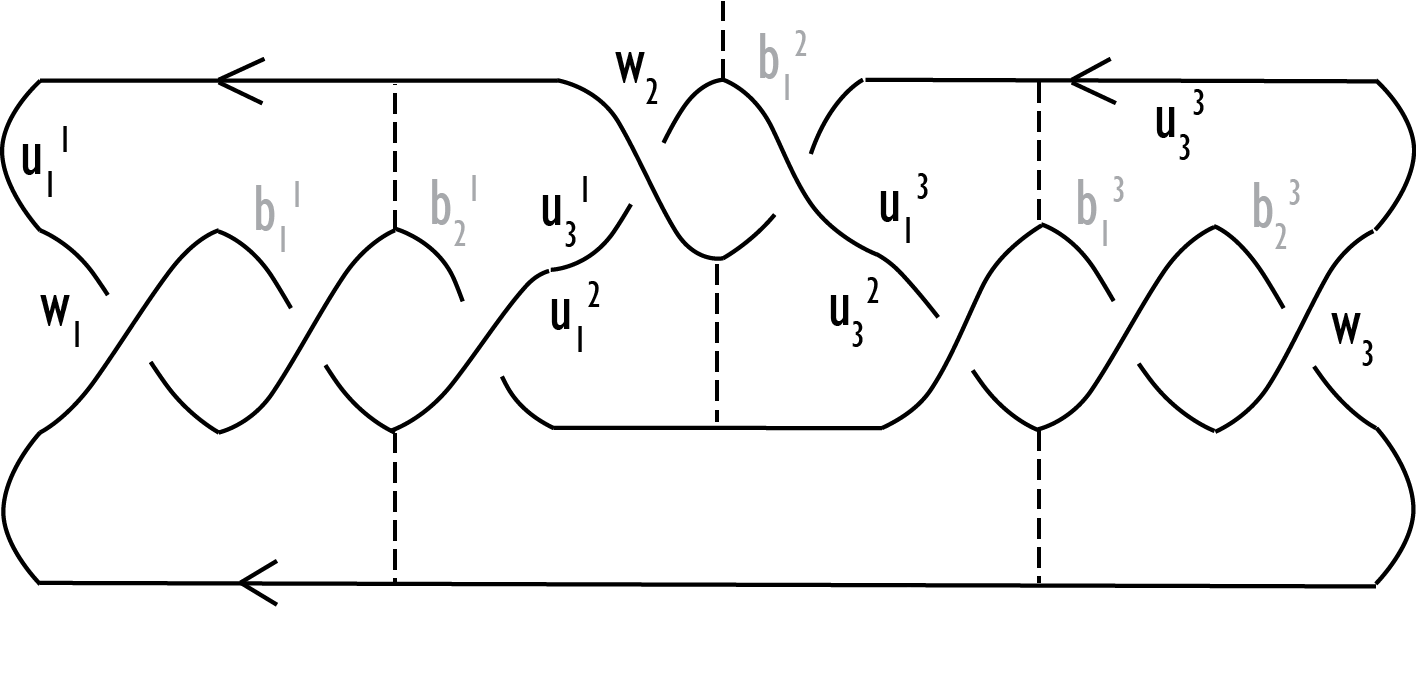}
\vspace{-0.15in}

\it Fig. 5
\end{center}

\ind It would be very interesting if one could generalize the process of obtaining the exact tetrahedra parameters from the diagram beyond
the family of two-bridged links. 

\vvv

\noind \textbf{4. Complex volume}
\vv

\ind The complex volume is an invariant of hyperbolic manifolds of mixed geometric and algebraic nature. It is a complex number: the real part is hyperbolic volume, and the imaginary part is the Chern-Simons invariant. (\cite{Neumann1992}).  

\v

Similar machinery and formulas from \cite{Zickert} can be used for a computation of the exact complex volume of a hyperbolic 2-bridge link \( L\). In particular, some of our edge labels (e.g., \( b_i^j\) ) correspond to labels of short edges (in the terminology of \cite{Zickert}) of tetrahedra in a triangulation of the complement of \( L\). To see the exact correspondence, faces and edges of every tetrahedron should be located on a reduced, alternating diagram of \( L\), and the gluing pattern should be traced. From these labels, one can compute a flattening of every tetrahedron, and then the complex volume of a complement of \( L\). Such a computation would save the step of developing an image of every cusp, allowing one to find the labels from a link diagram instead.

\vvv

\noind \textbf{5. Acknowledgments}
\vv

\ind  I would like to thank Makoto Sakuma for useful discussions, Christian Zickert for helpful correspondence, Oliver Dasbach for encouragement and interest in my work, and the referees, whose
comments helped to improve the paper. Particular thanks are due to Morwen Thistlethwaite for numerous enlightening conversations on the subject over time.

\vvvv
\noind Anastasiia Tsvietkova \newline
Mathematics Department \newline
University of California, Davis
One Shield Ave \newline
Davis, CA 95616, USA \newline
tsvietkova@math.ucdavis.edu


\newpage


\begin{thebibliography}{90}

\bibitem{Sakuma1} H. Akiyoshi, Hirotaka, M. Sakuma, M. Wada, Y. Yamashita, \textit{
Punctured torus groups and 2-bridge knot groups. I.}, 
Lecture Notes in Mathematics 1909, Springer, Berlin (2007), 252 pp.

\bibitem{Conway} J. Conway, \textit{An enumeration of knots and links, and some of their algebraic
properties}, Computational Problems in Abstract Algebra (Ed. Leech), Pergamon
Press (1967), 329--358.

\bibitem{Snap} D. Coulson, O. A. Goodman, C. D. Hodgson, W. D. Neumann, \textit{Computing arithmetic invariants of 3-manifolds}, Experiment. Math. 9 (2000), no. 1, 127--152.

\bibitem{Futer} D. Futer, E. Kalfagianni, J. S. Purcell, \textit{Dehn filling, volume, and the Jones polynomial}, J. Differential Geom. 78 (2008), no. 3, 429--464.

\bibitem{Futer2009} D. Futer, E. Kalfagianni, J. S. Purcell, \textit{Symmetric links and Conway sums: volume and Jones polynomial}, Math. Res. Lett. 16 (2009), no. 2, 233--253. 

\bibitem{Bounds} F. Gueritaud, \textit{On canonical triangulations of once-punctured torus bundles}, with
an appendix by D. Futer, Geom. Topol. 10 (2006), 1239--1284

\bibitem{Lackenby} M. Lackenby, \textit{The volume of hyperbolic alternating link complements}, with an appendix by I. Agol and D. Thurston., Proc. London Math. Soc. (3) 88 (2004), no. 1, 204--224. 

\bibitem{Neumann1992} W. D. Neumann, \textit{Combinatorics of triangulations and the Chern-Simons invariant for hyperbolic 3-manifolds}, Topology '90 (Columbus, OH, 1990), Ohio State Univ. Math. Res. Inst. Publ., 1, de Gruyter, Berlin (1992), 243--271.

\bibitem{WithWalter} W. D. Neumann, A. Tsvietkova, \textit{Intercusp geodesics and the invariant trace field of hyperbolic 3-manifolds}, in prep.

\bibitem{SakumaWeeks} M. Sakuma and J. R.
Weeks, \textit{Examples of canonical decomposition of hyperbolic link complements}, Japan. J. Math. (N. S.) 21 (1995), No. 2, 393--439.

\bibitem{ThistlethwaiteTsvietkova} M. Thistlethwaite, A. Tsvietkova, \textit{An alternative approach to hyperbolic structures on link complements}, to appear in Algebr. Geom. Topol.,  ArXiv: math.GT/1108.0510v1. 

\bibitem{Thurston} W. P. Thurston,
\textit{The Geometry and Topology of Three-Manifolds}, Electronic
Version 1.1 (March 2002), http://www.msri.org/publications/books/gt3m/  


\bibitem{Weeks}
J. R. Weeks, \textsl{SnapPea}: a computer program for
creating and studying hyperbolic $3$--manifolds, freely available
from


http://thames.northnet.org/weeks/index/SnapPea.html


 

\bibitem{Zickert} C. K. Zickert, \textit{The volume and Chern-Simons invariant of a representation}, Duke Math. J. 150 (2009), no. 3, 489--532.


\end{thebibliography}
\end{document}